# Fixed-time Cooperative Guidance Law for Multiple Missiles in Simultaneous Attacks

Hang Yu, Keren Dai, Haojie Li, Yao Zou*, Shaojie Ma and He Zhang*

**Abstract** To enable multiple missiles to attack a maneuvering target simultaneously, fixed-time distributed cooperative guidance laws are proposed in this paper. Here, we present a novel fixed-time fast nonsingular terminal sliding mode surface (FNTSMS). In particular, the sliding mode surface not only avoids singularities but also has the characteristic of a settling time boundary regardless of the initial conditions. Based on the FNTSMS, we have developed a distributed guidance law, which has the characteristic of fixed-time convergence. The guidance law achieves the consensus of range-to-go, relative velocity along and perpendicular to the line of sight (LOS) direction to realize the simultaneous attack. In addition, a saturation function is introduced to avoid the chattering problem caused by the commonly used sign function. Furthermore, the distributed cooperative guidance law with communication failure is considered and proved theoretically, which shows that the proposed guidance law has excellent performance. Finally, the simulation results verify the performance of the distributed guidance law and its robustness against communication topology mutations, and explain the phenomenon in detail.
Keywords: simultaneous attacks, distributed control, multiple missile.

## 1 Introduction

With the development of military technology, more and more important targets are equipped with anti-missile defense systems such as close-in weapon system (CIWS) [1, 2]. However, due to the existence of the defense system, for traditional single missile attack strategy, the defense system can intercept the attacker easily, which seriously affect the success rate and reliability rate of target attacks. Aiming at this, a strategy of multiple missiles attacking the target simultaneously was proposed[1].

According to the existed literature, there are two methods for multiple missiles simultaneous attack[3]:one is impact time control guidance(ITCG), set a common impact time for each missile before launched. And then let all missiles guide independently[1]. However, due to the different performance of each missile, it is difficult to complete the simultaneous attack.

The second method is to build a communication network so that each missile transmits its own information, such as time-to-go, range-to-go, velocity, heading angle and so on. By transmitting this information, each missile adjusts its own movement state through the guidance law, so that the missile states reach the consensus, and then simultaneous attack can be achieved. The established communication networks are divided into distributed networks and centralized networks. A Cooperative Proportional Navigation (CPN) was proposed[3], designed the time-varying navigation gain to reduce the variance of the time-to-go estimates of each missile, thereby the consensus of time-to-go estimates is achieved. The guidance law is still centralized and requires real-time updating of the missile information in the network.

Nevertheless, in a complex battlefield environment, the centralized network topology is difficult to guarantee. Therefore, in order to save missile energy and improve the missile's anti-jamming capability, it is necessary to develop distributed cooperative guidance laws. A distributed guidance law is presented that the time-to-go estimates consensus in finite-time[4]. But it omits that when the heading error is zero, the convergence time will tend to infinity. In [2] extended [4] and designed two distributed guidance laws to achieve the real time-to-go consensus in finite-time. In[5], a guidance law for finite-time convergence for static targets is proposed and considered that the consensus problem in directed networks under both fixed topology and switching topologies is investigated.

In order to realize a simultaneous attack, the agreement of the missile's times-to-go are usually selected as the basis for judging the consensus. However, since the relative velocity and distance between the missile and the target are changing in the complex battlefield environment, the missile's time-to-go of the missile cannot be obtained accurately. To overcome this difficulty, we choose the distance and velocity along the LOS direction as the consensus variable.

In [2-9], the guidance law is based on the convergence in finite-time, the settling time obtained by this guidance law depends on the initial conditions of the missile. However, in the actual battlefield environment, it is difficult to obtain accurate information about the initial conditions of the missile, thus, it is impossible to obtain an accurate settling time.

In view of the above, in [10], a control law of the convergence in fixed-time is proposed. It quickly attracted many scholars' research in multiple fields[11-17]. To the best knowledge of the authors, studies on the cooperative guidance law for simultaneous attack of multiple missiles and consensus in fixed-time are rare, except in [14].

In order to realize the fast response and robustness of multiple missile system, the sliding mode control is considered in the guidance law of this paper. A novel fast nonsingular terminal sliding mode surface(FNTSMS) is proposed[18, 19]. In this note, we proposed a distributed cooperative guidance law for Simultaneous attack. The guidance law based on FNTSMS, it provides the consensus of the range-to-go and heading error in fixed-time.

The key contributions in this paper include:

(1) A novel fast nonsingular terminal sliding mode surface (FNTSMS) is designed. The sliding mold surface not only avoids singularity, but also has the advantages of conventional sliding mold surface and the characteristics of convergence in fixed-time.

(2) The guidance law based on FNTSMS, it provides the consensus of the range-to-go and heading error in fixed-time, the boundary of settling time regardless of initial conditions. At the same time, introducing smooth function $\tanh(\cdot)$ in guidance law to reduce chatter.

(3) Discussing the impact of sudden changes in communication topology on

system stability during and after system convergence.

This paper is organized as follows. Section 2 states the preliminaries and engagement geometry. Distributed guidance law is mainly described in Section 3. Simulation analysis of the guidance law is presented in Section 4. Finally, conclusions are given in Section 5.

## 2 Preliminaries and engagement geometry

### 2.1 Preliminaries

The communication topology between missiles is described by the weighted graph $\mathcal{G}=(\mathcal{V}, \mathcal{E}, \mathcal{A})$。$\mathcal{V} = 1,2\ldots n$ is the set of nodes, an edge set $\varepsilon \subseteq \mathcal{V} \times \mathcal{V}$ and an adjacent matrix $\mathcal{A} = [a_{ij}] \in \mathbb{R}^{n\times n}$. If edge $(i, j) \in \varepsilon$, then node $i$ is called a neighbor of node $j$. An edge $(i, j)$ on $\mathcal{G}$ represents that the state of node $i$ is accessible to node $j$. The adjacent matrix $\mathcal{A}$ is defined as $a_{ij} = 1$ if edge $(j, i) \in \varepsilon$ and $a_{ij} = 0$ otherwise, especially $a_{ii} = 0$. The Laplacian matrix $\mathcal{L} = [l_{ij}] \in \mathbb{R}^{n\times n}$ is defined as $l_{ii} = \sum_{j=1}^{n} a_{ij}$ and $l_{ij} = -a_{ij}$ for $i \neq j$. A graph has a directed spanning tree if a subset of the edges forms a spanning tree.[2, 11].缺引用文献 [zuo]

The communication topology $\mathcal{G}$ satisfies the following assumption:

*Assumption* : $\mathcal{G}$ is connected, which means that, there exists a communication path that involves all the missiles.

*Lemma* 1 [20] Once for any $m \in \mathbb{R}^n$ satisfying $1^T m = 0$, we have $m^T \mathcal{L} m \geqslant \lambda_s m^T m$, where $\lambda_s$ denotes the smallest nonzero eigenvalue of $\mathcal{L}$.

*Lemma* 2 [11] Let $x_1, x_2, \ldots, x_n \geq 0$, and $0 < p_1 \leq 1, 1 < p_2 < \infty$. Then

$$\sum_{i=1}^{n} x_i^{p_1} \geq \left(\sum_{i=1}^{n} x_i\right)^{p_1} \tag{1}$$

$$\sum_{i=1}^{n} x_i^{p_2} \geq n^{1-p_2} \left(\sum_{i=1}^{n} x_i\right)^{p_2} \tag{2}$$

*Lemma* 3 [11] Consider the following system:

$$\dot{y} = -ay^{\frac{m}{n}} - by^{\frac{p}{q}}, \quad y(0)=y_0 \tag{3}$$

Where $a>0$, $b>0$, $m, n, p, q$ are positive odd integers and satisfying $m > n, p < q$. Then $y$ will converge to the origin in fixed-time, and the settling time bounded by:

$$T \le \frac{1}{a}\frac{n}{m-n} + \frac{1}{b}\frac{q}{q-p} \tag{4}$$

If $\vartheta \triangleq [q(m-n)/n(q-p)] \le 1$, then the conservative estimation of the settling time bounded by:

$$T \le \frac{q}{q-p}\left(\frac{1}{\sqrt{ab}}\tan^{-1}\sqrt{\frac{a}{b}} + \frac{1}{a\vartheta}\right) \tag{5}$$

*Definition* 1[10] For the multi-agent theory, if the settling time $T(X_0)$ is bounded, and the system equation is satisfied $\exists T_{max} > 0 : T(X_0) \le T_{max}, \forall X_0 \in \mathbb{R}^n$, then the system can be called globally fixed-time stable.

$$\lim_{t \to T}|x_i(t) - x_0(t)| \to 0$$

## 2.2 Engagement geometry

As shown in figure 1, it is the engagement geometry of *n* missiles attacks a target, where $V_i$ is the velocity of the *i*th missile and $r_i$ is the distance between the *i*th missile and the target. For a group of *n* missiles attack a target, once all missiles have already laterally headed on the target, the flight paths are only different in the vertical plane. In order to facilitate the study of the missile flight process, can be converted into a single missile attack on the target in planar engagement scenario.

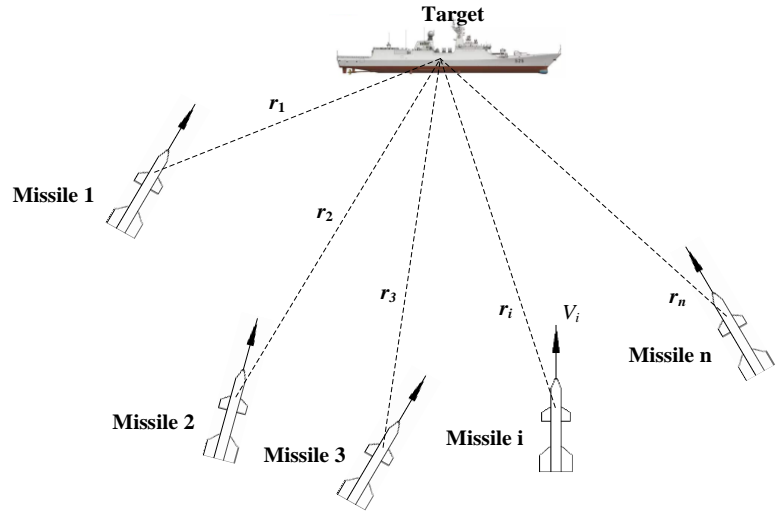

Fig.1 The engagement geometry of multiple missiles

According to the above description, the planar engagement scenario of the *i*th missile and the target is illustrated in Fig. 2. Subscript *i* is the variable related to the *i*th missile. The velocity of the *i*th missile along and perpendicular to the LOS direction are denoted by $V_{r,i}$, $V_{q,i}$, respectively. The normal acceleration, tangential acceleration, flight vehicle flight path angle, heading the, and LOS angle of the $i$th missile are denoted by $a_{n,i}$, $a_{t,i}$, $\lambda_i$, $\phi_i$, $\gamma_i$, respectively. Then the planar engagement kinematic model can be expressed as[2, 14]:

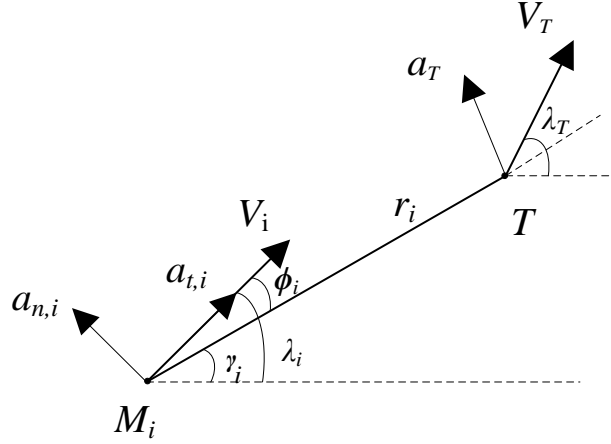

Fig. 2 The planar engagement scenario

$$\begin{aligned}
V_{r,i}(t) &= \dot{r}_i(t) \\
V_{q,i}(t) &= \dot{\lambda}_i(t) r_i(t) \\
\dot{V}_{r,i}(t) &= \frac{V_{q,i}^2(t)}{r_i(t)} - u_i^r(t) + u_{T,i}^r(t) \\
\dot{V}_{q,i}(t) &= -\frac{V_{q,i}(t) V_{r,i}(t)}{r_i(t)} - u_i^q(t) + u_{T,i}^q(t)
\end{aligned} \quad (6)$$

Where $u_i^q, u_i^r$ are the $i$th missile's acceleration components along and perpendicular to the LOS, respectively. For simplicity of expression, the time variable $t$ is omitted in subsequent articles.

## 3 Design of fixed-time cooperative guidance law against the maneuvering target

### 3.1 Cooperative Guidance Law Based on FNTSMS

The consensus errors can be expressed as:

$$\begin{cases} \xi_{r,i} = \sum_{j=1}^{N} a_{ij}\left(r_i - r_j\right) \\ \xi_{Vr,i} = \sum_{j=1}^{N} a_{ij}\left(V_{r,i} - V_{r,j}\right) \end{cases} \quad (7)$$

*Remark 1* For a group of missiles, once $\xi_{r,i}$ and $\xi_{Vr,i}$ are used to judge the consensus variable, if it satisfies $\xi_{r,i} = 0$, $\xi_{Vr,i} = 0$ for $t \geq T_{\max}$, then the group of missiles achieved the consensus. And the $T_{\max}$ is called the bound of settling time. Especially, compared with finite-time stability, the convergence time for fixed-time is independent of the initial state of the system.

In this subsection, the guidance law is designed with the help of FNTSMS. The main purpose of tangential acceleration $a_{t,i}$ is to adjust the missile's flight speed so that $r_i$ can achieve the consensus agreement. The purpose of normal acceleration $a_{n,i}$ is mainly to adjust the missile's flight direction, so that $V_{q,i}$ converges to zero, then to reduce the heading error $\phi_i$.

In this paper, a novel FNTSMS without any constraint is designed as follows:

$$s_i = \xi_{Vr,i} + \alpha_1 \xi_{r,i} + \alpha_2 \beta(\xi_{r,i}) \tag{8}$$

Where $\alpha_i > 0, i = 1, 2$, and $\beta(\xi_{r,i})$ is defined as:

$$\beta(\xi_{r,i}) = \begin{cases} -\dfrac{\alpha_1}{\alpha_2}\xi_{r,i} + \dfrac{\alpha_1}{\alpha_2}\xi_{r,i}^{m_1/n_1} + \xi_{r,i}^{p_1/q_1}, & \text{if } \bar{s}_i=0 \text{ or } \bar{s}_i \neq 0, |\xi_{r,i}| \geq \mu \\ l_1 \xi_{r,i} + l_2 \operatorname{sgn}(\xi_{r,i}) \xi_{r,i}^2, & \text{if } \bar{s}_i \neq 0, |\xi_{r,i}| < \mu \end{cases} \tag{9}$$

With $\bar{s}_i = \xi_{Vr,i} + \alpha_1 \xi_{r,i} + \alpha_2 \xi_{r,i}^{p_1/q_1}$, $l_1 = (2 - p_1/q_1)\mu^{p_1/q_1 - 1}$, $l_2 = (p_1/q_1 - 1)\mu^{p_1/q_1 - 2}$,

$m_1/n_1 > 1$, $p_1/q_1 < 1$, $\mu > 0$ is a small constant.

With exchanging the information between neighboring missiles via a communication network, a novel distributed guidance law is proposed as:

$$\begin{aligned} a_{n,i} &= u_i^q \cos\phi_i - u_i^r \sin\phi_i \\ a_{t,i} &= u_i^r \cos\phi_i + u_i^q \sin\phi_i \end{aligned} \tag{10}$$

Where $u_i^r$ and $u_i^q$ given by:

$$u_i^r = \frac{1}{\sum_{j=1}^n a_{ij}} \left\{ \sum_{j=1}^n a_{ij}\left(\frac{V_{q,i}^2}{r_i} - \frac{V_{q,j}^2}{r_j} + u_j^r\right) + (\alpha_1 \xi_{Vr,i} + \alpha_2 \dot{\beta}(\xi_{r,i})) + k_1 s_i^{\frac{m_2}{n_2}} + k_2 s_i^{\frac{p_2}{q_2}} + \eta_1 \tanh(s_i) \right\} \tag{11}$$

$$u_i^q = -\frac{V_{q,i} V_{r,i}}{r_i} + k_3 V_{q,i}^{\frac{m_1}{n_1}} + k_4 V_{q,i}^{\frac{p_1}{q_1}} + \varepsilon |V_{q,i}|^\omega \operatorname{sign}(V_{q,i}) \tag{12}$$

Where $k_i > 0, i = 1,2,3,4$, $m_2/n_2 > 1$, $p_2/q_2 < 1$, $\eta_1 > 0$, $\varepsilon > 0$, $1 > \omega > 0$.

*Theorem* 1 Consider the multiple missile system (6) with the distributed guidance law (10), then $s_i = 0 (i=1,2,\ldots,n)$ can be achieved in fixed-time and settling time is derived

as:

$$T_1 \leq \frac{n^{(m_2-n_2)/2n_2}}{k_1} \frac{n_2}{m_2 - n_2} + \frac{1}{k_2} \frac{q_2}{q_2 - p_2} \tag{13}$$

*Proof* Differentiating $s_i$ in (8) against time yields

$$\dot{s}_i = \dot{\xi}_{Vr,i} + \alpha_1 \xi_{Vr,i} + \alpha_2 \beta(\xi_{r,i}) \tag{14}$$

$$\dot{\beta}(\xi_{r,i}) = \begin{cases} (-\dfrac{\alpha_1}{\alpha_2} + \dfrac{m_1 \alpha_1}{n_1 \alpha_2} \xi_{r,i}^{m_1/n_1 - 1} + \dfrac{p_1}{q_1} \xi_{r,i}^{p_1/q_1 - 1}) \xi_{Vr,i}, & \text{if } \bar{s}_i = 0 \text{ or } \bar{s}_i \neq 0, |\xi_{r,i}| \geq \mu \\ (l_1 + 2l_2 \operatorname{sgn}(\xi_{r,i}) \xi_{r,i}) \xi_{Vr,i}, & \text{if } \bar{s}_i \neq 0, |\xi_{r,i}| < \mu \end{cases} \tag{15}$$

Differentiating $\xi_{r,i}$ and $\xi_{Vr,i}$ in (7) against time, then substituting the result into (15), it yields

$$\begin{aligned}\dot{s}_i &= \sum_{j=1}^{n} a_{ij} \left(\frac{V_{q,i}^2}{r_i} - \frac{V_{q,j}^2}{r_j} - u_i^r + u_j^r\right) + \alpha_1 \xi_{Vr,i} + \alpha_2 \dot{\beta}(\xi_{r,i}) \\ &= -\sum_{j=1}^{n} a_{ij} u_i^r + \sum_{j=1}^{n} a_{ij} \left(\frac{V_{q,i}^2}{r_i} - \frac{V_{q,j}^2}{r_j} + u_j^r\right) + \alpha_1 \xi_{Vr,i} + \alpha_2 \dot{\beta}(\xi_{r,i})\end{aligned} \tag{16}$$

Substituting (11) into above result, then equation (16) can be rewritten as

$$\dot{s}_i = \sum_{i}^{n} \left(-k_1 s_i^{\frac{m_2}{n_2}} - k_2 s_i^{\frac{p_2}{q_2}} - \eta_1 \tanh(s_i)\right)$$

Consider the following Lyapunov function

$$V_1 = \frac{1}{2} \sum_{i=1}^{n} s_i^2 \tag{17}$$

Differentiating (17) against time, then substituting $\dot{s}_i$ into it, it yields

$$\begin{aligned}\dot{V}_1 &= \sum_{i=1}^{n} s_i \dot{s}_i \\ &= \sum_{i=1}^{n} s_i \left(-k_1 s_i^{\frac{m_2}{n_2}} - k_2 s_i^{\frac{p_2}{q_2}} - \eta_1 \tanh(s_i)\right)\end{aligned} \tag{18}$$

According to the properties of $\tanh(\cdot)$ function, there is the following relation

$$s_i \tanh(s_i) = s_i \frac{e^{s_i} - e^{-s_i}}{e^{s_i} + e^{-s_i}} = \frac{s_i(e^{2s_i} - 1)}{e^{2s_i} + 1} \tag{19}$$

and

$$\begin{cases} s_i(e^{2s_i} - 1) \geq 0, s_i \geq 0 \\ s_i(e^{2s_i} - 1) < 0, s_i < 0 \end{cases}$$

Consider the result of equation (19), the equation (18) becomes

$$\dot{V}_1 \le \sum_{i=1}^{n} -k_1 (s_i^2)^{\frac{m_2+n_2}{2n_2}} - k_2 (s_i^2)^{\frac{p_2+q_2}{2q_2}} \tag{20}$$

According to Lemma 2, inequality (20) becomes

$$\dot{V}_1 \le -k_1 n^{\frac{n_2-m_2}{2n_2}} \left( \sum_{i=1}^{n} s_i^2 \right)^{\frac{m_2+n_2}{2n_2}} - k_2 \left( \sum_{i=1}^{n} s_i^2 \right)^{\frac{p_2+q_2}{2q_2}} \tag{21}$$

And $V_1 = (1/2)\sum_{i=1}^{n} s_i^2$, substituting it into the preceding inequality obtains

$$\dot{V}_1 \le -k_1 n^{\frac{n_2-m_2}{2n_2}} (2V_1)^{\frac{m_2+n_2}{2n_2}} - k_2 (2V_1)^{\frac{p_2+q_2}{2q_2}} \tag{22}$$

Let $y = \sqrt{2V_1}$, taking it into, then

$$\dot{y} \le -k_1 n^{\frac{n_2-m_2}{2n_2}} y^{\frac{m_2}{n_2}} - k_2 y^{\frac{p_2}{q_2}} \tag{23}$$

According to Lemma 3, it can be obtained that $s_i$ will converge to 0 in fixed-time, and the bound of convergence time is

$$T_1 \le \frac{n^{(m_2-n_2)/2n_2}}{k_1} \frac{n_2}{m_2-n_2} + \frac{1}{k_2} \frac{q_2}{q_2-p_2} \tag{24}$$

From the above proof, we can get that when $s_i = 0$, $i=1,2,\ldots,n$, for $t \ge T_1$. Therefore, when $t \ge T_1$, the multiple missiles system will enter the sliding mode and move towards the origin. By selecting appropriate parameters $k_1, k_2, m_2, n_2, p_2, q_2$ the convergence rate and control effort cost are optimized.

*Remark* 2 According to Lemma 3, once the $\vartheta_2 \triangleq [q_2(m_2-n_2)/n_2(q_2-p_2)] \le 1$, then the conservative estimation of the settling time bounded by:

$$T_1 \le \frac{q_2}{q_2-p_2} \left( \frac{1}{\sqrt{k_1 n^{(n_2-m_2)/2n_2} k_2}} \tan^{-1} \sqrt{\frac{k_1 n^{(n_2-m_2)/2n_2}}{k_2}} + \frac{1}{k_1 n^{(n_2-m_2)/2n_2} \vartheta_2} \right) \tag{25}$$

After the proof of Theorem 1, we can get $s_i = 0$ in fixed-time, but in order to make the multiple missiles system converge to the origin, that is $\xi_{r,i} = 0$ and $\xi_{Vr,i} = 0$, need to prove next.

*Theorem* 2 On the basis of Theorem 1's results $s_i = 0$, according to the sliding mode surface, then $\xi_{r,i} = 0$ and $\xi_{Vr,i} = 0$ will be achieve in fixed time $T_2$. The settling time is derived as

$$T_2 \le T_1 + \frac{n^{(m_1-n_1)/2n_1}}{\alpha_1} \frac{n_1}{m_1-n_1} + \frac{1}{\alpha_2} \frac{q_1}{q_1-p_1} \tag{26}$$

*Proof* Substituting $s_i = 0$ and the first equation of (9) into (8), it yields

$$\dot{\xi}_{Vr,i} = -(\alpha_1 \xi_{r,i}^{m_1/n_1} + \alpha_2 \xi_{r,i}^{p_1/q_1}) \tag{27}$$

Consider the following Lyapunov function

$$V_2 = \frac{1}{2} \sum_{i=1}^{n} \xi_{r,i}^2 \tag{28}$$

Differentiating (28) against time, it yields

$$\dot{V}_2 = \sum_{i=1}^{n} \xi_{r,i} \dot{\xi}_{Vr,i} \tag{29}$$

Then substituting (29) into above result and according to the lemma 2, becomes

$$\begin{aligned}
\dot{V}_2 &= -\sum_{i=1}^{n} \xi_{r,i} \left( \alpha_1 \xi_{r,i}^{\frac{m_1}{n_1}} + \alpha_2 \xi_{r,i}^{\frac{p_1}{q_1}} \right) \\
&= -\alpha_1 \sum_{i=1}^{n} \xi_{r,i}^{\frac{m_1+n_1}{n_1}} - \alpha_2 \sum_{i}^{n} \xi_{r,i}^{\frac{p_1+q_1}{q_1}} \\
&\leq -\alpha_1 n^{\frac{n_1-m_1}{2n_1}} \left( \sum_{i=1}^{n} \xi_{r,i}^2 \right)^{\frac{m_1+n_1}{2n_1}} - \alpha_2 \left( \sum_{i=1}^{n} \xi_{r,i}^2 \right)^{\frac{p_1+q_1}{2q_1}}
\end{aligned} \tag{30}$$

Substituting $V_2 = (1/2)\sum_{i=1}^{n} \xi_{r,i}^2$ into the last inequality of

$$\dot{V}_2 \leq -\alpha_1 n^{\frac{n_1-m_1}{2n_1}} (2V_2)^{\frac{m_1+n_1}{2n_1}} - \alpha_2 (2V_2)^{\frac{p_1+q_1}{2q_1}} \tag{31}$$

Let $y = \sqrt{2V_1}$, taking it into

$$\dot{y} \leq -\alpha_1 n^{\frac{n_1-m_1}{2n_1}} y^{\frac{m_1}{n_1}} - \alpha_2 y^{\frac{p_1}{q_1}} \tag{32}$$

According to Lemma 3, it can be obtained that $\xi_{r,i}$ will converge to the origin in fixed-time, and the bound of convergence time is

$$T_2 \leq T_1 + \frac{n^{(m_1-n_1)/2n_1}}{\alpha_1} \frac{n_1}{m_1 - n_1} + \frac{1}{\alpha_2} \frac{q_1}{q_1 - p_1} \tag{33}$$

*Remark* 3 Same as Remark 2, when $\vartheta_1 \triangleq [q_1(m_1 - n_1)/n_1(q_1 - p_1)] \leq 1$, the bound of convergence time can be obtained instead as

$$T_2 \leq T_1 + \frac{q_1}{q_1 - p_1} \left( \frac{1}{\sqrt{\alpha_1 \alpha_2 n^{(n_1-m_1)/2n_1}}} \tan^{-1} \sqrt{\frac{\alpha_1 n^{(n_1-m_1)/2n_1}}{\alpha_2}} + \frac{1}{\alpha_1 n^{(n_1-m_1)/2n_1} \vartheta_1} \right) \tag{34}$$

*Remark* 4 The sliding mode surface $s_i=0$ has been proven in theorem 1, and $\xi_{r,i} = 0$ has been proven in theorem 2. Thus $\xi_{Vr,i} = 0$ for $t \geq T_2$. And the bound of convergence time $T_1$ and $T_2$ are independent of the initial state of the system. This is a very unique feature compared to finite-time convergence.

*Theorem* 3 Consider the multiple missile system (6) with the distributed guidance law (13), then $V_{q,i}=0$ ($i$=1,2,…,n) can be achieved in fixed-time and settling time is derived as:

$$T_3 \leq \frac{n^{(m_3-n_3)/2n_3}}{k_3} \frac{n_3}{m_3-n_3} + \frac{1}{k_3} \frac{q_3}{q_3-p_3} \quad (35)$$

*Proof* Substituting the (13) into last equation of (6), it yields

$$\dot{V}_{q,i}(t) = -k_3 V_{q,i}^{\frac{m_3}{n_3}} - k_4 V_{q,i}^{\frac{p_3}{q_3}} - \varepsilon |V_{q,i}|^{\omega} \operatorname{sign}(V_{q,i}) \quad (36)$$

Choose the following Lyapunov function

$$V_3 = \frac{1}{2}\sum_{i=1}^{n} V_{q,i}^{2} \quad (37)$$

Differentiating (37) against time, and then Substituting (36) into it, becomes

$$\begin{aligned}
\dot{V}_3 &= \sum_{i=1}^{n} \left( -k_3 V_{q,i}^{\frac{m_3+n_3}{n_3}} - k_4 V_{q,i}^{\frac{p_3+q_3}{q_3}} - \varepsilon |V_{q,i}|^{\omega} V_{q,i} \operatorname{sign}(V_{q,i}) \right) \\
&= -k_3 \sum_{i=1}^{n} V_{q,i}^{\frac{m_3+n_3}{n_3}} - k_4 \sum_{i=1}^{n} V_{q,i}^{\frac{p_3+q_3}{q_3}} - \varepsilon |V_{q,i}|^{\omega+1} \\
&\leq -k_3 \sum_{i=1}^{n} V_{q,i}^{\frac{m_3+n_3}{n_3}} - k_4 \sum_{i=1}^{n} V_{q,i}^{\frac{p_3+q_3}{q_3}}
\end{aligned} \quad (38)$$

The rest of the proof is similar to that of theorem 1 or theorem 2, which is omitted for brevity.

*Remark* 5 After the above proof, $V_{q,i}$ will converge to the origin in fixed-time, which means that the heading error $\phi_i$ will also converge to the origin in fixed-time. Therefore, the missile will accurately attack the target.

## 4 Simulation examples

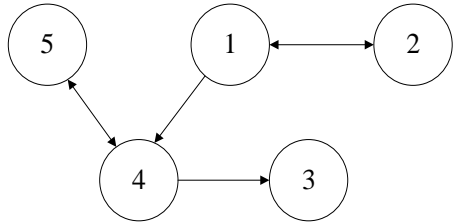

Fig. 3 The information flow $\mathcal{G}$

In this section, simulation results are reported to illustrate the effectiveness of the proposed guidance law (10) and (43). Moreover, the robustness of the communication topology is considered during and after the convergence of the multi-missile system. Consider a network of five missiles communicating with each other according to Fig. 3. The network of five missiles is connected and undirected. The LOS of distance, velocities of the missiles and the heading error are listed in Table 1.

In the simulations, we consider three different communication topologies, as shown in Fig. 3

The simulation is divided into two parts. First, the graph shown in Figure 3 is used to verify the performance of the proposed guidance law 10. Second, the three different communication topologies shown in Figure 4 are used to verify the robustness of the guidance law 10 communication topology.

Table 1 Initial states of each missile

| missile | LOS of distance ($r$, m) | Initial speed ($V$, m/s) | heading error ($\Phi$, rad) |
| --- | --- | --- | --- |
| missile 1 | 16000 | 350 | -0.09 |
| missile 2 | 15050 | 320 | 0.10 |
| missile 3 | 13990 | 270 | 0.11 |
| missile 4 | 13950 | 300 | -0.15 |
| missile 5 | 15000 | 331 | 0.12 |

The parameters of the guidance laws (10) and the FNTSMS (8) are given as:

$k_1 = 0.265$, $k_2 = 2$, $k_3 = 0.25$, $k_4 = 2$, $p_2 = 5$, $q_2 = 7$, $m_2 = 13$, $n_2 = 11$, $\varepsilon = 1.5$, $\omega = 0.5$, $\eta_1 = 1.5$, $\mu = 0.001$, $p_1 = 5$, $q_1 = 7$, $m_1 = 11$, $n_1 = 9$, $\alpha_1 = 0.25$, $\alpha_2 = 2$, $p_3 = 3$, $q_3 = 5$, $m_3 = 7$, $n_3 = 5$

The parameters of the guidance laws (43) are given as:

$k_1 = 2$, $k_2 = 3$, $k_3 = 0.1$, $p_3 = 3$, $q_3 = 5$, $m_3 = 7$, $n_3 = 5$, $N_0 = 4$, $\eta_2 = 1$

The simulation is divided into two parts. First, the guidance law (10), then the guidance law (43) is simulated.

**4.1 Simulation of guidance law (10)**

As shown in Fig. 4, the simulation results of the guidance law (10). It is clear that the guidance law can effectively complete a simultaneous attack on a stationary target. It can be found from Fig. 4(a), the time taken for each missile to arrive the target is 51.0s. As revealed from Fig. 4(b) and Fig. 4(c), the consensus error of ranges to go and the consensus error of velocities in the direction of LOS converges to zero is 19.9s and 21.2s, respectively. Since $\vartheta_1 < 1$ and $\vartheta_2 < 1$, then according to (25) and (34), the settling time bounded by $T_2 \leq 48.21$s, the convergence time shown in Figure 4(b) meets the boundary conditions. According to the (35), we can get the settling time $T_3 \leq 14.9$s, it can be verified in Fig. 4(d), in this picture, time of the consensus error of the heading errors converges to zero is 2.8s. According to the above results, the designed guidance law (10) can meet the requirements.

The initial position, velocity and acceleration components of the target are (15000,15000), $V_x = 100, V_y = 100, a_x = 3.5\sin(0.5t + \frac{7}{6}\pi)$, $a_y = 3.5\sin(0.5t + \frac{7}{6}\pi)$, respectively. which are time-varying and unknown to the missiles.

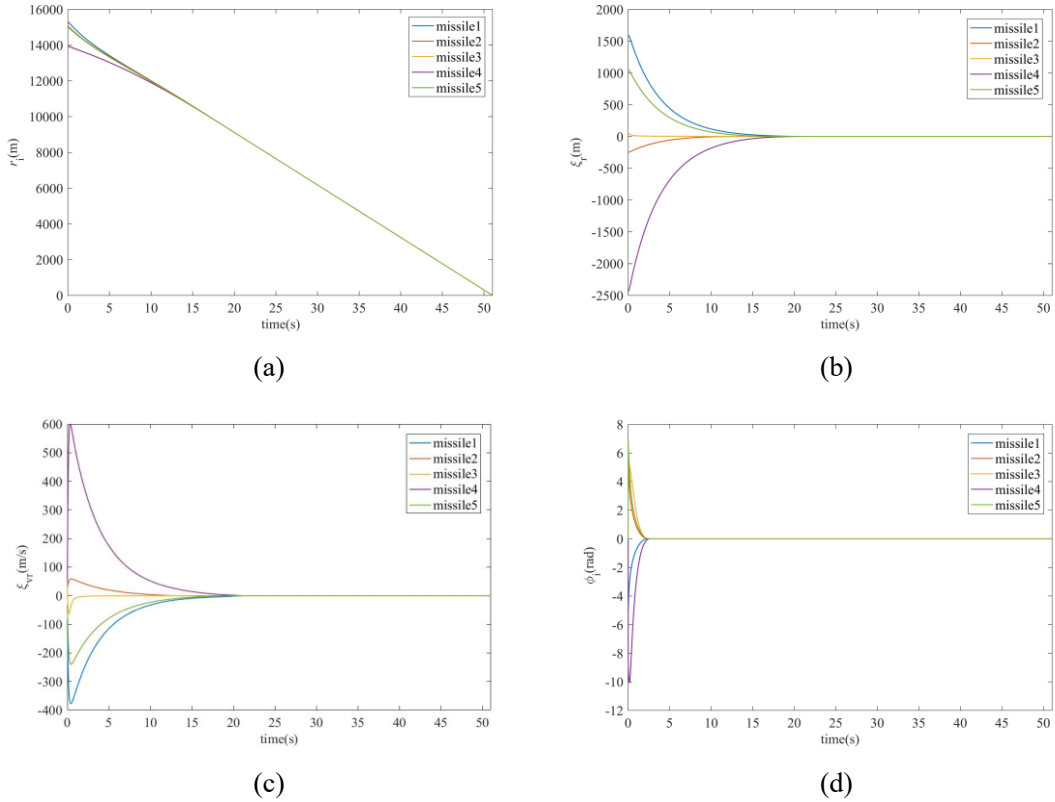

Fig. 4 The guidance law (10): (a) Missiles range-to-go; (b) The consensus error of ranges to go; (c) The consensus error of velocities along LOS; (d) The heading errors;

Compared with guidance law (10), guidance law (43) has smaller convergence time than before, so is the arrival time. The simulation results are consistent with the theoretical analysis.

## 5 Conclusion

This article proposes that a fixed-time convergence distributed cooperative guidance law for attacking a stationary target. Here, a novel fast nonsingular terminal sliding mode surface (FNTSMS) without any constraint is designed. Particularly, the sliding mold surface not only avoids singularity, but also has the advantages of traditional sliding mold surface and the characteristics of convergence in fixed-time. Based on the proposed FNTSMS, we have developed the guidance law, which has the characteristics of fixed-time convergence. Due to the introduction of the tanh(·) function, system chatter is avoided. Further on, it is verified by simulation that the guidance law 10 can perform a cooperative attack on stationary and maneuvering targets, and the robustness of the communication topology is verified.

# References:

[1]. Impact-Time-Control Guidance Law for Anti-Ship Missiles.

[2]. Jialing, Z. and Y. Jianying, Distributed Guidance Law Design for Cooperative Simultaneous Attacks with Multiple Missiles.

[3]. Homing Guidance Law for Cooperative Attack of Multiple Missiles.

[4]. Simultaneous attack of a stationary target using multiple missiles: a consensus-based approach.

[5]. Consensus of leader-followers system of multi-missile with time-delays and switching topologies.

[6]. Cooperative Guidance Law Design for Simultaneous Attack with Multiple Missiles Against a Maneuvering Target.

[7]. Finite time convergence cooperative guidance law based on graph theory.

[8]. Finite-time cooperative guidance laws for multiple missiles with acceleration saturation constraints.

[9]. Integrated guidance and control law for cooperative attack of multiple missiles.

[10]. Andrey, P., Nonlinear Feedback Design for Fixed-Time Stabilization of Linear Control Systems. IEEE TRANSACTIONS ON AUTOMATIC CONTROL, 2012.

[11]. Zuo, Z., Nonsingular fixed-time consensus tracking for second-order multi-agent networks. Automatica.

[12]. Alireza Khanzadeh, M.P., Fixed-time sliding mode controller design for synchronization of complex dynamical networks. 2017.

[13]. Wang, C., et al., Fixed-Time Formation Control of Multirobot Systems: Design and Experiments. IEEE Transactions on Industrial Electronics, 2019. 66(8): p. 6292-6301.

[14]. Guofei Li, Y.W.P.X., Fixed-time Cooperative Guidance Law with Input Delay for Simultaneous Arrival. International Journal of Control.

[15]. Zuo, Z., Q. Han and B. Ning, An Explicit Estimate for the Upper Bound of the Settling Time in Fixed-Time Leader-Following Consensus of High-Order Multivariable Multiagent Systems. IEEE Transactions on Industrial Electronics, 2019. 66(8): p. 6250-6259.

[16]. A New Class of Finite-Time Nonlinear Consensus Protocols for Multi-Agent Systems. 2014.

[17]. Distributed robust finite-time nonlinear consensus protocols for multi-agent systems. 2014.

[18]. Kunfeng Lu, Y.X., Adaptive attitude tracking control for rigid spacecraft with finite-time convergence. 2013.

[19]. Olfati-Saber, R., Consensus problems in networks of agents with switching topology and time-delays.